\documentclass[11pt]{amsart}
\usepackage{amssymb}
\usepackage{amscd}
\usepackage{amsfonts}
\usepackage{amsmath}
\usepackage{verbatim}

\theoremstyle{plain}
\newtheorem{definition}{Definition}[section]
\newtheorem{lemma}[definition]{Lemma}
\newtheorem{thm}[definition]{Theorem}
\newtheorem{prop}[definition]{Proposition}
\newtheorem{cor}[definition]{Corollary}

\makeindex
\begin{document}
\title{Maximally homogeneous para-CR manifolds of semisimple type}
\author{D.V. Alekseevsky, C. Medori\and A. Tomassini}
\date{}
\address{Department of Mathematics\\ Hull University\\
UK}
\email{d.v.alekseevsky@hull.ac.uk}
\address{Dipartimento di Matematica\\ Universit\`a di Parma\\ Viale G.\,P. Usberti, 53/A\\ 43100
  Parma\\ Italy}
\email{costantino.medori@unipr.it}
\address{Dipartimento di Matematica\\ Universit\`a di Parma\\ Viale G.\,P. Usberti, 53/A\\ 43100
  Parma\\ Italy}
\email{adriano.tomassini@unipr.it}
\thanks{This work was supported by the Project M.U.R.S.T. ``Geometric Properties of
Real and Complex Manifolds'' and by G.N.S.A.G.A.
of I.N.d.A.M. The first author was also supported by Grant FWI
Project P17108-N04 (Wien) and the Leverhulme Trust, EM/9/2005/0069.}
\begin{abstract}
An almost para-CR structure on a
manifold $M$ is given by a distribution $HM \subset TM$ together with a
field $K \in \Gamma({\rm End}(HM))$ of involutive endomorphisms of
$HM$. If $K$ satisfies an integrability condition, then $(HM,K)$ is
called a para-CR structure. The notion of maximally homogeneous
para-CR structure of a semisimple type is given. A
classification of such maximally homogeneous para-CR
structures is given in terms of appropriate gradations of real
semisimple Lie algebras.
\end{abstract}
\maketitle
%%%%%%%%%%%%%%%%
\tableofcontents
%%%%%%%%%%%%%%%%% SIMBOLI MATEMATICI %%%%%%%%%%%%%%%%%%%%%%%%%%%%%%%%%%%%%%%%%%%%%%%%%%

%%%%%%%%%%%%%%% GOTICHE%%%%%%%%%%%%%%%%%%%%%%%

\def\gA{{\mathfrak A}} \def\ga{{\mathfrak a}}  \def\cA{{\mathcal A}}
\def\gB{{\mathfrak B}} \def\gb{{\mathfrak b}}  \def\cB{{\mathcal B}}
\def\gC{{\mathfrak C}} \def\gc{{\mathfrak c}}  \def\cC{{\mathcal C}}
\def\gD{{\mathfrak D}} \def\gd{{\mathfrak d}}  \def\cD{{\mathcal D}}
\def\gE{{\mathfrak E}} \def\gge{{\mathfrak e}}  \def\cE{{\mathcal E}}
\def\gF{{\mathfrak F}} \def\gf{{\mathfrak f}}  \def\cF{{\mathcal F}}
\def\gG{{\mathfrak G}} \def\ggg{{\mathfrak g}} \def\cG{{\mathcal G}}
\def\gH{{\mathfrak H}} \def\gh{{\mathfrak h}}  \def\cH{{\mathcal H}}
\def\gI{{\mathfrak I}} \def\gi{{\mathfrak i}}  \def\cI{{\mathcal I}}
\def\gJ{{\mathfrak J}} \def\gj{{\mathfrak j}}  \def\cJ{{\mathcal J}}
\def\gK{{\mathfrak K}} \def\gk{{\mathfrak k}}  \def\cK{{\mathcal K}}
\def\gL{{\mathfrak L}} \def\gl{{\mathfrak l}}  \def\cL{{\mathcal L}}
\def\gM{{\mathfrak M}} \def\gm{{\mathfrak m}}  \def\cM{{\mathcal M}}
\def\gN{{\mathfrak N}} \def\gn{{\mathfrak n}}  \def\cN{{\mathcal N}}
\def\gO{{\mathfrak O}} \def\go{{\mathfrak o}}  \def\cO{{\mathcal O}}
\def\gP{{\mathfrak P}} \def\gp{{\mathfrak p}}  \def\cP{{\mathcal P}}
\def\gQ{{\mathfrak Q}} \def\gq{{\mathfrak q}}  \def\cQ{{\mathcal Q}}
\def\gR{{\mathfrak R}} \def\gr{{\mathfrak r}}  \def\cR{{\mathcal R}}
\def\gS{{\mathfrak S}} \def\gs{{\mathfrak s}}  \def\cS{{\mathcal S}}
\def\gT{{\mathfrak T}} \def\gt{{\mathfrak t}}  \def\cT{{\mathcal T}}
\def\gU{{\mathfrak U}} \def\gu{{\mathfrak u}}  \def\cU{{\mathcal U}}
\def\gV{{\mathfrak V}} \def\gv{{\mathfrak v}}  \def\cV{{\mathcal V}}
\def\gW{{\mathfrak W}} \def\gw{{\mathfrak w}}  \def\cW{{\mathcal W}}
\def\gX{{\mathfrak X}} \def\gx{{\mathfrak x}}  \def\cX{{\mathcal X}}
\def\gY{{\mathfrak Y}} \def\gy{{\mathfrak y}}  \def\cY{{\mathcal Y}}
\def\gZ{{\mathfrak Z}} \def\gz{{\mathfrak z}}  \def\cZ{{\mathcal Z}}

\def\gli{{\mathfrak {gl}}}
\def\sli{{\mathfrak {sl}}}
\def\so{{\mathfrak {so}}}
\def\su{{\mathfrak {su}}}
\def\sp{{\mathfrak {sp}}}

%%%%%%%%%%%%%%%%%%%% GRUPPI %%%%%%%%%%%%%%%%%%%%%%%%%%%%%
\newcommand{\Id}{{\rm I}}
\newcommand{\GL}{{\bf GL}}
\newcommand{\SO}{{\bf SO}}
\newcommand{\SU}{{\bf SU}}
\newcommand{\Gra}{{\bf G \!\bf r}}%%% compare 2volte nel cap5(grassmann)
\newcommand{\rad}{{\hbox{rad}}}
\newcommand{\Rad}{{\mathcal R}} %nella sezione3 def. il sist. di radici
\newcommand{\tr}{{\hbox{tr}}}

%%%%%%%%%%%%%%%%%% R,C,Z %%%%%%%%%%%%%%%%%%%%%%%

\newcommand\C{{\mathbb C}}
\newcommand\R{{\mathbb R}}
\newcommand\Z{{\mathbb Z}}
\newcommand\T{{\mathbb T}}

\newcommand\pCR{{para-CR\,\,}}
\newcommand\bC{{\bf C}}

%%%%%%%%%%%%%%%%%%%%%%%%%%%%%%%%%%%%%%%%%%%%%%%%%%%%%%
%\section*{Introduction}
\section{Introduction and notation} Let $M$ be a $2n$-dimensional manifold. An
\emph{almost paracomplex structure}\index{paracomplex structure} on $M$ is a field of
endomorphisms  $K \in {\rm End}(TM)$ of the tangent bundle $TM$ of
$M$ such  that $K^2=\hbox{\rm id}$. It is called an (almost) paracomplex
structure in the {\em strong sense} if  its $\pm 1$-eigenspace
distributions
$$
T^\pm M=\{X\pm KX\,\,\vert\,\,X\in\Gamma(M,TM)\}
$$
have the same rank (see e.g. \cite{L}, \cite{CFG}). An almost
paracomplex structure $K$ is called a \emph{paracomplex
structure}\index{paracomplex structure}, if it is
\emph{integrable}, i.e.
$$
S(X,Y)=[X,Y]+[KX,KY] - K[X,KY] - K[KX,Y]=0\
$$
for any vector fields $X,\,Y\in \Gamma(TM)$.\newline
This is equivalent to say that the distributions
$T^\pm M$ are involutive.\smallskip\par
Recall that an \emph{almost}
CR-\emph{structure}\index{CR-structure} of codimension $k$ on a $2n+k$-dimensional
manifold $M$ is a distribution $HM\subset TM$ of rank $2n$
together with a field of endomorphisms $J \in {\rm End}(HM)$ such
that $J^2=-\hbox{\rm id}$.
An almost CR-structure is called
CR-\emph{structure}\index{CR-structure}, if the $\pm i$-eigenspace subdistributions
$H_\pm M$ of the complexified tangent bundle $T^\C M$ are
involutive.\newline We define an almost para-CR structure in a
similar way.
\begin{definition}
An \emph{almost para-CR structure}\index{para-CR structure}
of codimension $k$ on a $2n+k$-dimensional
manifold $M$ (in the weak sense) is a pair $(HM,K)$, where $HM \subset TM$ is a rank
$2n$ distribution and $ K \in {\rm End}(HM)$ is a field of
endomorphisms such that $K^2=\hbox{\rm id}$ and $K\neq\pm\hbox{\rm
id}$.\newline
An  almost para-CR structure is said to be a
{\rm para-CR structure}\index{para-CR structure}, if the
eigenspace subdistributions $H_\pm M \subset HM$
are integrable or equivalently if  the following integrability conditions hold:
\begin{eqnarray}
&{}&[KX,KY]+[X,Y]\in\Gamma(HM)\,,\label{integrability1}\\[3pt]
&{}& S(X,Y):=[X,Y]+[KX,KY] -
K([X,KY]+[KX,Y])=0\label{integrability2}
\end{eqnarray}
for all  $X,\,Y\in\Gamma(HM)$.
\end{definition}
If the eigenspace distributions
$$
H_\pm M =\{X\pm KX\,\,\vert\,\, X\in\Gamma(M,HM)\}
$$
of an almost para-CR structure have the same rank, then $(HM,K)$ is called an
almost para-CR structure in the {\em strong sense}.
\newline
A straightforward computation shows that the
integrability condition is equivalent to the involutivness of the
distributions $H_+M$ and $H_-M$.
\newline A manifold $M$, endowed with an (almost) para-CR
structure, is called an (\emph{almost})
\emph{para-{\rm CR} manifold}\index{para-CR manifold}
\index{para-CR manifold}.\newline
Note that a direct product of (almost) para-CR manifolds is
an (almost) para-CR manifold.\newline
One can associate with a  point $x \in M$ of a para-CR manifold $(M,HM,K)$
a fundamental graded Lie algebra $\gm$. A para-CR structure is said to be {\em regular} if these Lie algebras
$\gm_x$ do not depend on $x$. In this case, a para-CR structure can be considered as a Tanaka structure
(see \cite{AS} and section 4). A regular para-CR structure is called a structure of {\em semisimple type}
if  the full prolongation
$$
\ggg = \gm^{\infty} = \gm^{-d} + \cdots + \gm^{-1} + \ggg^0 + \ggg^1 + \cdots
$$
  of  the
 associated non-positively graded Lie algebra $\ggg^{-d} + \cdots + \ggg^{-1} + \ggg^0$
  (which is an analogue of the generalized Levi form
 of a CR structure) is a semisimple Lie algebra. Such a para-CR
 structure defines a parabolic geometry and its group of
 automorphisms ${\rm Aut}(M,HM,K)$ is a Lie group of dimension $\leq \dim
 \ggg$.\newline
 Recently in \cite{NS} P. Nurowski and G. A. J. Sparling consider the natural para-CR structure which arises
 on the $3$-dimensional space $M$ of solutions of a second order ordinary differential
 equation $y''=Q(x,y,y')$. Using the Cartan method of prolongation, they construct the full prolongation
 $\mathcal G\to M$ of $M$ with a $\sli(3,\R)$-valued Cartan connection and a quotient line bundle over $M$
 with a conformal metric of signature $(2,2)$. This is a para-analogue of the Feffermann bundle of a CR-structure. They
 apply these bundles to the initial ODE and get interesting applications. \newline
 In \cite{AMT} we proved that a para-CR structures of semisimple type
on a simply connected manifold $M$ with  the automorphism group of maximal dimension $\dim \ggg$
can be identified with a (real) flag manifold $M = G/P$ where $G$ is the simply connected Lie group
with the Lie algebra $\ggg$  and $P$ the parabolic subgroup
generated by the parabolic subalgebra $\gp = \ggg^0 + \ggg^1 + \cdots +\ggg^d.$
    We gave a classification of maximally homogeneous para-CR
 structures of semisimple type such that the associated graded semisimple Lie algebra $\ggg$
has depth $d=2$. In the present paper we classify all maximally homogeneous para-CR structures of
semisimple type in terms of graded real semisimple Lie algebras.

%%%%%%%%%%%%%%%%%%%%%%%%%%%%%%%%%%%%%%%%%%%%%%%%%%%%%%%%%%%%%%%%%%%
\section{Graded Lie algebras associated with para-CR structures}
\subsection{Gradations of a Lie algebra} Recall that a \emph{gradation}\index{gradation}
(more precisely a $\Z$-gradation) of
\emph{depth} $k$ of a Lie algebra $\ggg$ is a direct sum decomposition
\begin{equation}\label{gradation}
\ggg =\sum_{i\in\Z}\ggg^{i}= \ggg^{-k}+\ggg^{-k+1}+\cdots+\ggg^{0}+
\cdots+\ggg^{j}+\cdots
\end{equation}
such that $[\ggg^{i},\ggg^{j}]\subset \ggg^{i+j}$,
for any $i,j\in \Z$, and $\ggg^{-k}\neq \{0\}$. Note that $\ggg^0$ is a
subalgebra of $\ggg$ and each $\ggg^i$ is a $\ggg^0$-module.\newline
We say that an element $x\in\ggg^{j}$
has \emph{degree} $j$ and we write $d(x)=j$. The
endomorphism $\delta$ of $\ggg$ defined by
$$
\delta_{\vert_{\ggg_{j}}}=j\cdot id
$$
is a  semisimple derivation of $\ggg$ (with integer eigenvalues), whose
eigenspaces determine the gradation. Conversely, any semisimple
derivation  $\delta$ of $\ggg$ with integer eigenvalues defines a
gradation where the grading space $\ggg^j$ is the eigenspace of
$\delta$ with eigenvalue $j$. If $\ggg$ is a semisimple Lie
algebra, then any derivation $\delta$ is inner, i.e. there exists $d\in\ggg$ such that
$\delta=\hbox{\rm ad}_d$. The element $d\in\ggg$ is called the \emph{grading element}\index{grading element}.
\begin{definition} A gradation $\ggg = \sum \ggg^i$ of a Lie
algebra {\rm (}or a graded Lie algebra $\ggg${\rm )}
is called
\begin{enumerate}
\item     {\rm fundamental}\index{graded Lie algebra!fundamental}, if the negative part
$\gm = \sum_{i< 0} \ggg^i $ is
generated by $\ggg^{-1}$;
 \item  {\rm (almost) effective} or {\rm transitive}\index{graded Lie algebra!transitive},
if the non-negative part
$$
\ggg^{\geq 0}=\gp=\ggg^0+\ggg^1+\cdots
$$
contains no non-trivial ideal of $\ggg$;
 \item { \rm non-degenerate}\index{graded Lie algebra!non-degenerate}, if
$$
X\in\ggg^{-1}\,,\,\,[X,\ggg^{-1}]=0\,\,\Longrightarrow\,\,X=0\,.
$$
\end{enumerate}
\end{definition}
\subsection{Fundamental algebra associated with a distribution}
Let ${\mathcal H}$ be a distribution on a manifold $M$. We recall that to any point $x\in M$ it is possible
to associate a Lie algebra $\gm(x)$ in the following way.\newline
First of all, we  consider  a filtration of the Lie algebra ${\mathcal X}(M)$ of
vector fields defined inductively by
\begin{eqnarray*}
\Gamma({\mathcal H})_{-1}&=&\Gamma({\mathcal H})\,,\\
\Gamma({\mathcal H})_{-i}&=&\Gamma({\mathcal H})_{-i+1}+
[\Gamma({\mathcal H}),\Gamma({\mathcal H})_{-i+1}]\,\,,\hbox{\rm for}\,\,\,i>1.
\end{eqnarray*}
Then evaluating vector fields at a point $x\in M$, we get a flag
$$
T_xM\supset\mathcal{H}_{-d-1}(x)=\mathcal{H}_{-d}(x)
\supsetneq \mathcal{H}_{-d+1}(x)
\supset \cdots
\supset {\mathcal H}_{-2}(x)\supset \mathcal{H}_{-1}(x)=\mathcal{H}_x
$$
in $T_xM$, where
$$
{\mathcal
H}_{-i}(x)=\{X_{\vert_x}\,\,\vert\,\,X\in\Gamma({\mathcal
H})_{-i}\}.
$$
Let us assume that $\mathcal{H}_{-d}(x)=T_xM$. The commutators of vector fields induce a structure of
fundamental negatively graded Lie algebra on the associated graded
space
$$
{\gm}(x)={\rm gr}(T_xM)={\gm}^{-d}(x)+{\gm}^{-d+1}(x)+\cdots + {\gm}^{-1}(x)\,,
$$
where ${\gm}^{-j}(x)={\mathcal H}_{-j}(x)/{\mathcal H}_{-j+1}(x)$. Note that
${\gm}^{-1}(x)={\mathcal H}_x$. \newline
A distribution ${\mathcal H}$ is called a \emph{regular distribution}\index{regular distribution} of
\emph{depth} $d$ and \emph{type}
$\gm$ if all graded Lie algebras ${\gm}(x)$ are isomorphic to a given
graded fundamental Lie algebra
$$
{\gm}={\gm}^{-d}+ {\gm}^{-d+1}+\cdots + {\gm}^{-1}\,.
$$
In this case ${\gm}$ is called the \emph{Lie algebra associated} with the
 distribution ${\mathcal H}$.
A regular distribution $\mathcal{H}$ is called \emph{non-degenerate}
if the associated  Lie algebra is non-degenerate.
 %%%%%%%%%%%%%%%%%%%%%%%%%%%%%%%%%%%%%%%%%%%%%%
%%%%%%%%%%%%%%%%%%%%%%%%%%%%%%%%%%%%%%%%%%%%%%%%%%%%%%%%%%%%%%%%%%
%%%%%%%%%%%%%%%%%%%%%%%%%%%%%%%%%%%%%%%%%%%%%%%%%%%%%%%%%%%%%%%%%%
\subsection{Para-CR algebras and regular para-CR structures} We recall the
following
\begin{definition} A pair $(\gm,K_o)$, where $\gm = \gm^{-d} + \cdots +
\gm^{-1}$ is a negatively graded fundamental Lie algebra and $K_o$
is an involutive endomorphism of $\gm^{-1}$, is called a
\emph{para-{\rm CR} algebra}\index{para-CR algebra}
of \emph{depth} $d$. If, moreover, the $\pm
1$-eigenspaces $\gm^{-1}_{\pm}$ of $K_o$ on $\gm^{-1}$ are
commutative subalgebras of $\gm$, then $(\gm,K_o)$ is called an
{\em integrable para-CR algebra}\index{para-CR algebra!integrable}.
\end{definition}
\begin{definition}\label{regularstructure}
Let $(\gm, K_o)$   be a para-CR algebra of depth $d$.
An almost  para-CR structure $(HM,K)$
on $M$ is called \emph{regular} of
\emph{type} $(\gm, K_o)$ and \emph{depth} $d$ if, for any $x\in
M$, the pair $({\gm}(x),K_x)$ is isomorphic to $({\gm},K_o)$. We
say that the regular almost para-CR structure is non-degenerate
if the graded algebra $\gm$ is non-degenerate.
\end{definition}
Note that a regular almost para-CR structure of type $(\gm ,K_0)$ is integrable
if and only if the Lie algebra $(\gm ,K_0)$ is integrable.
%%%%%%%%%%%%%%%%%%%%%%%%%%%%%%%%%%%%%%%%%%%%%%%%%%%%%%%%%%%%%%%%%%%%%%%%
\section{Prolongations of graded Lie algebras}%%%%%%%%%%%%%%%%%%%%%%%%%%%%%%%%
\subsection{Prolongations of negatively graded Lie algebras}%%%%%%%%%%%%%%%
The \emph{full prolongation}\index{full prolongation} of a negatively graded fundamental Lie algebra
$\gm=\gm^{-d}+\cdots + \gm^{-1}$ is defined as a maximal
graded Lie algebra
$$
\ggg(\gm)=\ggg^{-d}(\gm)+\cdots +\ggg^{-1}(\gm)+
\ggg^0(\gm)+\ggg^1(\gm)+\cdots
$$
with the negative part
$$
\ggg^{-d}(\gm)+\cdots +\ggg^{-1}(\gm)=\gm
$$
such that the following transitivity condition holds: \newline
$$
\mbox{if}\,\, X\in\ggg^k(\gm)\,,\,k\geq 0\,,\,\,
[X,\ggg^{-1}(\gm)]=\{0\}\,,\,\,\mbox{then}\,\,\, X=0\,.
$$
In \cite{T1}, N. Tanaka proved that the full prolongation $\ggg(\gm)$
always exists and it is unique up to isomorphisms. Moreover, it can be
defined inductively by
$$
\ggg^i(\gm)=
\begin{cases}
\gm^i & \hbox{\rm if}\,\, i<0\,,\\
\{A\in\hbox{\rm Der}(\gm,\gm)\,:\, A(\gm^j)\subset \gm^j\,,\forall j<0\}
& \hbox{\rm if}\,\, i=0\,,\\
\{A\in\hbox{\rm Der}(\gm,\sum_{h< i}\ggg^h(\gm))\,:\,
A(\gm^j)\subset \ggg(\gm)^{i+j}\,,\forall j<0\}
& \hbox{\rm if}\,\, i>0\,,
\end{cases}
$$
where ${\rm Der}(\gm, V)$ denotes the space of derivations of the
Lie algebra $\gm$ with values in the $\gm$-module $V$.

Note that
\begin{eqnarray}
\ggg^i(\gm)&\!\!=\!\!&
\Big\{A\in\hbox{\rm Hom}_\R(\gm,\sum_{h<i}\ggg^h(\gm))\,\Big\vert\,
A(\ggg^h(\gm))\subset\ggg^{h+i}(\gm)\,\,\forall h<0\,,\,\,\,\\
\nonumber &{}&\mbox{and}\,\, [A(Y),Z]+ [Y,A(Z)]=A([Y,Z])\,\,\forall Y,Z\in\gm
\Big\}\,.
\end{eqnarray}
%%%%%%%%%%%%%%%%%%%%%%%%%%%%%%%%%%%%%%%%%%%%%%%%%%%%%%%%%%%%%%%%
%%%%%%%%%%%%%%%%%%%%%%%%%%%%%%%%%%%%%%%%%%%%%%%%%%%%%%%%%%%%%%%%
%%%%%%%%%%%%%%%%%%%%%%%%%%%%%%%%%%%%%%%%%%%%%%%%%%%%%%%%%%%%%%%%
\subsection{Prolongations of non-positively graded Lie algebras}
Consider now a non-positively graded Lie algebra $\gm
+\ggg^0=\gm^{-d}+\cdots +\gm^{-1}+\ggg^0$. The \emph{full
prolongation} of $\gm +\ggg^0$ is the subalgebra
$$
(\gm+\ggg^0)^\infty
=\gm^{-d}+\cdots +\gm^{-1} +\ggg^0+\ggg^1+\ggg^2+\cdots
$$
of $\ggg(\gm)$, defined inductively by
$$
\ggg^i=\{X\in\ggg(\gm)^i\,:\,[X,\gm^{-1}]\subset\ggg^{i-1}\}\,,\,\,\,
\hbox{\rm for any}\,\, i\geq 1\,.
$$
\begin{definition}
A graded Lie algebra ${\gm}+{\ggg}^0$ is called of {\rm finite
type}\index{graded Lie algebra!of finite type} if its full prolongation
$\ggg=({\gm}+{\ggg}^0)^\infty$ is a
finite dimensional Lie algebra and it is called of {\rm semisimple
type}\index{graded Lie algebra!of semisimple type} if $\ggg$ is a
finite dimensional semisimple Lie algebra.
\end{definition}
We have the following criterion (see \cite{T2}, \cite{AS})
\begin{lemma}\label{criterion}
Let $(\gm =\sum_{i<0}\gm^{i}, K_o)$ be an integrable para-CR
algebra and $\ggg^0$ the subalgebras of $\ggg^0(\gm)$ consisting
of any $A\in\ggg^0(\gm)$ such that $A\vert_{\gm^{-1}}$ commutes
with $K_o$. Then the graded Lie algebra  $(\gm + \ggg^0)$ is of
finite type if and only if $\gm$ is non-degenerate.
\end{lemma}
The following result will be used in the last section (see e.g. \cite{MN1}, Theorem 3.21)
\begin{lemma}
Let $\ggg =\sum_i \ggg_i$ be a fundamental effective semisimple
graded Lie algebra such that $\gm +\ggg^0$ is of finite type. Then
$\ggg$ coincides  with the full prolongation $(\gm +\ggg^0)^\infty$
of $\gm +\ggg^0$.
\end{lemma}
%\begin{lemma} \label{non-degenerate lemma} A fundamental effective graded
%semisimple Lie algebra $\ggg = \sum \ggg^{i}$
%is  degenerate if and only if it contains a graded (not necessary
%proper) ideal of depth one.
%\end{lemma}
%%%%%%%%%%%%%%%%%%%%%%%%%%%%%%%%%%%%%%%%%%%%%%%%%%%%%
%%%%%%%%%%%%%%%%%%%%%%%%%%%%%%%%%%%%%%%%%%%%%%%%%%%%%
\section{Standard almost para-CR manifolds}
\subsection{Maximally homogeneous Tanaka structures}A regular para-CR
structure of type $(\gm ,K_0)$ is of
{\em finite type}\index{para-CR structure!of finite type} or,
respectively, of {\em semisimple type}\index{para-CR structure!of semisimple type},
if the Lie algebra $(\gm
+ \ggg^0)^{\infty}$ is finite-dimensional or, respectively,
semisimple. Recall that $\ggg^0 = Der(\gm, K_0)$ is the Lie
algebra of
Lie group ${\rm Aut}(\gm,K_0)$. \\
We recall the following (see \cite{AS})
\begin{definition} Let $\gm = \gm^{-d}+ \cdots + \gm^{-1} $ be a
negatively graded Lie algebra generated by $\gm^{-1}$ and $G^0$
a closed Lie subgroup of (grading preserving) automorphisms of
$\gm$. A \emph{Tanaka structure}\index{Tanaka structure} of \emph{type}
 $(\gm, G^0)$ on a manifold $M$ is a regular distribution
$\mathcal H \subset TM$ of type $\gm$ together with a principal
$G^0$-bundle $\pi : Q \to M$ of adapted coframes of $\mathcal H$.
A coframe $\varphi : {\mathcal H}_x \to \gm^{-1}$ is called
\emph{adapted} if it can be extended to an isomorphism $\varphi :
\gm_x \to  \gm$ of Lie algebra.
\end{definition}
We say that the Tanaka structure of type $(\gm , G^0)$ is of
\emph{finite type}\index{Tanaka structure!of finite type}
(respectively \emph{semisimple type}\index{Tanaka structure!of semisimple type}
$(\gm, G^0)$),
if the graded Lie algebra ${\gm}+\ggg^0$ is of finite type
(respectively semisimple type). Let $P$ be a  Lie subgroup of a
connected Lie group $G$ and $\gp$
(respectively $\ggg$) the Lie algebra of $P$ (respectively $G$).
\begin{thm}
Let $(\pi : Q\to M,{\mathcal H})$ be a Tanaka structure on $M$ of
semisimple type $(\gm,G^0)$. Then the Tanaka prolongation of
$(\pi,{\mathcal H})$ is a $P$-principal bundle ${\mathcal G}\to M$,
with the parabolic structure group $P$, equipped with a Cartan
connection $\kappa: T{\mathcal G}\to \ggg$, where $\ggg$ is the full
prolongation of $\gm +\ggg^0$ and ${\rm Lie}P = \gp =\sum_{i\geq 0}
\ggg_i$. Moreover, $\hbox{\rm Aut} ({\mathcal H},\pi)$ is a Lie
group and
$$
\dim\hbox{\rm Aut} ({\mathcal H},\pi)\leq \dim \ggg\,.
$$
\end{thm}
Let $({\mathcal H}, \pi:Q\to M)$ be a Tanaka structure of
semisimple type $(\gm,G^0)$ and $\ggg=(\gm+\ggg^0)^\infty=\gm
+\gp$ be the full prolongation of the non-positively graded Lie
algebra $\gm+\ggg^0$.
 \begin{definition} A semisimple  Tanaka structure $({\mathcal H}, \pi:Q\to M)$
is  called  {\rm maximally homogeneous}\index{maximally homogeneous} if the dimension of its
automorphism
group $ {\rm Aut}({\mathcal H}, \pi)$ is equal to  $\dim \ggg$.\\
\end{definition}
\subsection{Tanaka structures of semisimple type}
We construct maximally
homogeneous Tanaka structures of semisimple type $(\gm, G^0)$ as
follows. Let $G=\hbox{\rm Aut}(\ggg)$ be the Lie group of
automorphisms of the graded Lie algebra $\ggg$. Recall that $G^0$ is a
closed subgroup of the automorphism group of the graded Lie algebra
$\ggg^-=\gm$. Since the Lie algebra $\ggg$ is canonically associated
with $\gm$, we can canonically extend the action of  $G^0$ on $\gm$
to the action of $G^0$ on $\ggg$ by automorphisms. In other
words, we have an embedding $G^0\hookrightarrow \hbox{\rm
Aut}(\ggg)=G$ as a closed subgroup. We denote by $G^+$ the connected
(closed) subgroup of $G$ with Lie algebra $\ggg_+=\sum_{p>0}\ggg^p$.
Then $P=G^0\cdot G^+ \subset G$ is a (closed) parabolic subgroup of
$G$. Let $G/P$ be the corresponding flag manifold. We have a
decomposition $\ggg =\gm +\gp$ and we identify $\gm$ with the
tangent space $T_o(G/P)$. Then the natural action of $G^0$ on $\gm$ is
the isotropy representation of $G^0$. We have a natural Tanaka
structure $({\mathcal H},\pi :Q\to G/P)$ of type $(\gm,G^0)$,
where ${\mathcal H}$ is the $G$-invariant distribution defined by
$\gm^{-1}$ and $Q$ is the $G^0$-bundle of coframes on ${\mathcal H}$.\newline
\begin{comment}
which is generated by coframes
$$
Q_{\vert_o}=G^0\cdot \hbox{\rm id}\,,\quad \hbox{id}:\gm^{-1}\to \gm^{-1}\,.
$$
\end{comment}
 Hence, the flag manifold $G/P$ carries a natural maximally
 homogeneous Tanaka structure $(\mathcal H, \pi: Q \to G/P)$.\\
 The universal covering $F$ of the manifold $G/P$  also has the
 induced
 Tanaka structure  $(H_F, \pi_F: Q_F \to F)$ of type $(\gm,G^0 )$
 and the simply connected (connected) Lie group $\tilde G$ with
 the Lie algebra $\ggg$ acts  transitively and almost effectively on $F$ as a group
 of automorphisms of this Tanaka structure. Moreover, the
 stabilizer in  $\tilde G$ of an appropriate point $o \in F$ is the
 (connected)  parabolic  subgroup $\tilde P$
 generated by the subalgebra $\gp = \ggg^0 + \ggg^1 + \cdots +
 \ggg^d$.
  The Tanaka structure $(\mathcal{H}, \pi: Q \to F = \tilde G/\tilde P)$ is
obviously maximally homogeneous and it is called the
\emph{standard (simply connected maximally homogeneous) Tanaka structure} of type
$(\gm ,G^0)$. We can state the following (see e.g. \cite[Theor. 4.8]{AMT})
\begin{thm} \label{maxhomogeneous thm}
 Any maximally homogeneous Tanaka structure of  semisimple type
 $(\mathfrak{m}, G_0)$ is isomorphic to the standard Tanaka
 structure on the simply connected flag manifold $F = \tilde G / \tilde P$
 where $\tilde G$ is the simply connected semisimple Lie group
 with the Lie algebra $\ggg = (\gm + \ggg^0)^\infty$ and $\tilde P$
 is the parabolic subgroup generated by the subalgebra
 $\gp = \ggg^0 + \ggg^1 + \cdots +
 \ggg^d$.
\end{thm}

%\begin{thm} \cite{AMT} \label{maxhomogeneous thm}  There exists unique
%(up to an isomorphism)$  simply connected maximally homogeneous
%Tanaka structure of a given semisimple type  (\mathcal{H},
%\pi: Q \to F = G/P)$ be the standard maximally homogeneous Tanaka
%structure  of a semisimple type $(\gm, G^0)$. Any maximally
%homogeneous  Tanaka structure $(\mathcal{H}', \pi' : Q' \to M =
%G'/P')$   of a semisimple  type   is locally
%isomorphic to the standard maximally homogeneous Tanaka structure
%$(\mathcal{H}, \pi: Q \to F = G/P)$.
%\end{thm}

   Let $(HM,K)$ be a regular almost para-CR structure of
type $(\gm,K_0)$. Assume that it has finite type, i.e. $m=\dim
(\gm + \ggg^0)^{\infty} < \infty$. According to the above
definition,  $(HM,K)$ is { \it maximally homogeneous},
if it admits a (transitive) Lie group of automorphisms of dimension $m$.\\
By Theorem \ref{maxhomogeneous thm}, a maximally homogeneous
almost para-CR structure of semisimple type is locally equivalent
to the standard structure associated with a gradation of a semisimple
Lie algebra. In the following subsection we describe this
correspondence in more details.
\subsection{Models of almost para-CR manifolds}
%%%%%%%%%%%%%%%%%%%%%%
%%%%%%%%%%%%%%%%%%%%%%%%%%%%%%%%%%%%%%%%%%%%%%%%%%%%%%%%%%%%%%%%%%%%%%%%%%%%%%
%%%%%%%%%%%%%%%%%%%%%%%%%%%%%%%%%%%%%%%%%%%%%%%%%%%%%%%%%%%%%%%%%%%%%%%%%%%%%%
Let $\ggg = \sum_{-d}^{d} \ggg^i = \ggg^{-} + \ggg^0 + \ggg^+$ be
an effective fundamental  gradation of a semisimple Lie algebra
$\ggg$ with negative part $\gm = \ggg^{-}=\sum_{i<0}\ggg^i$ and positive part
$\ggg^+=\sum_{i>0}\ggg^i$.\newline\noindent Denote by $F = \tilde G/\tilde P$
the simply connected  real flag
manifold associated with the graded Lie algebra $\ggg$ where $\tilde G$
is the simply connected  Lie group with Lie algebra $\ggg$ and
$\tilde P = G^0 G^+$ is the connected subgroup generated by the
Lie subalgebra $\ggg^0 + \ggg^+$.\newline\noindent
We will identify
the tangent space $T_oF$ at the point $o=eP$ with the subspace
$$
\ggg/\gp\simeq\gm\,.
$$
Since the subspace $(\ggg^{-1}+\gp)/\gp\subset T_oF$ is invariant
under the isotropy representation of $P$, it defines an invariant
distribution $\mathcal{H}$ on $F$. Since the gradation is
fundamental, one can easily check that, for any $x\in F$, the
negatively graded Lie algebra $\gm(x)$ associated with ${\mathcal H}$
is isomorphic to the Lie algebra $\gm$. Moreover, let
\begin{equation}\label{decomposition(-1)}
    \ggg^{-1} = \ggg^{-1}_+ + \ggg^{-1}_-
\end{equation}
be a decomposition of  the
$G^0$-module  $\ggg^{-1}$ into a sum of two submodules and $K_0$
the associated  ${\rm ad}_{\ggg_0}$-invariant endomorphism such that
$\ggg^{-1}_{\pm}$ are $\pm 1$-eigenspaces of $K_0$.\newline
The decomposition (\ref{decomposition(-1)}) defines two
invariant complementary subdistributions
${\mathcal H}_\pm$ of the distribution ${\mathcal H}\subset TF$
associated with $\ggg^{-1}$ and $K_0$ defines $\tilde G$-invariant
para-CR structure  $(HF, K)$ on $F$. It is the
standard para-CR structure associated with the graded Lie
algebra  $\ggg$ and the decomposition (\ref{decomposition(-1)}). We
get the following theorem (see also \cite[Theor. 5.1]{AMT})
\begin{thm}\label{reductiontheorem}
Let $F=\tilde G/\tilde P$ be the simply connected flag manifold associated
with a (real) semisimple
effective fundamental graded Lie algebra $\ggg$. A decomposition
$\ggg^{-1}=\ggg^{-1}_+ +\ggg^{-1}_-$ of $\ggg^{-1}$ into
complementary $G^0$-submodules $\ggg^{-1}_\pm$ determines an
invariant almost para-CR structure $(HM,K)$ such that $\pm1
$-eigenspaces $H_\pm M $ of $K$ are subdistributions of $ HM$
associated with $\ggg^{-1}_\pm$. Conversely, any standard
almost para-CR structure $( HM, K)$ on $F$ can be obtained in such a way.\\
Moreover, $(HM,K)$ is:
\begin{enumerate}
    \item an almost para-CR structure if $\ggg^{-1}_+$
and $\ggg^{-1}_-$ have the same dimensions,
    \item  a para-CR
structure if and only if $\ggg^{-1}_+$ and $\ggg^{-1}_-$ are
commutative subalgebras of $\ggg$,
    \item non-degenerate if and only if $\ggg$ has no graded ideals
    of depth one.
\end{enumerate}
\end{thm}
%\begin{definition}\label{admissiblegradation} A fundamental effective  gradation of a semisimple Lie
%$algebra $\ggg$ is called {\rm admissible } if the $\ggg^0$-module
%$\ggg^{-1}$ is reducible.
%\end{definition}

  By Theorem \ref{reductiontheorem},
 the classification of maximally homogeneous para-CR
structures of semisimple type, up to local isomorphisms (i.e. up
to coverings), reduces to the description of all
gradation of semisimple  Lie algebras $\ggg$ and to decomposition of
the $\ggg^0$-module $\ggg^{-1}$ into irreducible submodules. We
will give such a description for complex and real semisimple Lie
algebras in the next two sections.
%%%%%%%%%%%%%%%%%%%%%%%%%%%%%%%%%%%%%%%%%%%%%%%%%%%%%%%%%%%%%%%%%%%%%%%%%%%%%%%%%%%%%%%%%%%%
%%%%%%%%%%%%%%%%%%%%%%%%%%%%%%%%%%%%%%%%%%%%%%%%%%%%%%%%%%%%%%%%%%%%%%%%%%%%%%%%%%%%%%%%%%%%%
%%%%%%%%%%%%%%%%%%%%%%%%%%%%%%%%%%%%%%%%%%%%%%%%%%%%%%%%%%%%%%%%%%%%%%%%%%%%%%%%%%%%%%%%%%%%
\section{Fundamental gradations of a complex semisimple Lie algebra} We recall here the
construction of a gradation of a complex semisimple Lie algebra $\ggg$.
Let $\gh$ be a Cartan subalgebra of a semisimple Lie algebra
$\ggg$
and
$$
\ggg =\gh \oplus\sum_{\alpha\in R}\ggg_\alpha
$$
be the root decomposition of $\ggg$ with respect to $\gh$. We denote by
$$
\Pi=\{\alpha_1,\ldots ,\alpha_\ell\}\subset R
$$
a system of simple roots of the root system $R$ and
associate to each simple root $\alpha_i$ (or corresponding vertex
of the Dynkin diagram) a non-negative integer $d_i$.
Using the {\em label vector}\index{label vector} $\vec{d}=(d_1,\ldots , d_\ell)$, we define the
\emph{degree} of a root $\alpha =\sum_{i=1}^\ell k_i\alpha_i$ by
$$
d(\alpha)=\sum_{i=1}^\ell k_id_i\,.
$$
This defines a gradation of $\ggg$ by the conditions
$$
d(\gh)=0\,,\qquad
d(\ggg_\alpha)=d(\alpha),\quad \forall\alpha\in R\,,
$$
which is called the \emph{gradation associated with the label
vector} $\vec{d}$. \newline
We denote by $d\in\gh$ the
corresponding grading element. Then $d(\alpha)=\alpha(d)$. Any
gradation of a complex semisimple Lie algebra $\ggg$ is conjugated
to a gradation of such a type (see \cite{GOV}). In particular, it
has the form
$$
\ggg=\ggg^{-k}+\cdots +\ggg^0+\cdots + \ggg^{k}\,,
$$
where $\ggg^0$ is a reductive subalgebra of $\ggg$ and the grading
spaces $\ggg^{-i}$ and $\ggg^{i}$ are dual with respect to the
Killing form. It is clear now that any graded semisimple Lie
algebra is a direct sum of graded simple Lie algebras. Hence, it
is sufficient to
describe gradations of simple Lie algebras.\\
We need the following (see \cite{Y})
\begin{lemma}
The gradation of a complex semisimple Lie algebra $\ggg$
associated with a label vector $\vec{d}=(d_1,\ldots , d_\ell)$ is fundamental
if and only if all labels $d_i\in\{0,1\}$.
\end{lemma}
Let $\Pi^1\subset\Pi$ be a set of simple roots. We denote by
$\vec{d}_{\Pi^1}$ the label vector which associates label one
to the roots in $\Pi^1$ and label zero to the other simple
roots.\newline Now we describe the depth of a fundamental
gradation.\newline
Let $\mu$ be the maximal root with respect to
the fundamental system $\Pi$. It can be written as a linear
combination
\begin{equation}\label{maximalroot}
\mu = m_1\alpha_1 +\cdots + m_\ell\alpha_\ell
\end{equation}
of fundamental roots, where the coefficient $m_i$ is a positive
integer called the
\emph{Dynkin mark associated} with $\alpha_i$.
\begin{lemma}\label{depth}
Let $\Pi^1=\{\alpha_{i_1},\ldots , \alpha_{i_s}\}\subset\Pi$ be a
set of simple roots.
Then the depth $k$ of the fundamental gradation defined by the label vector
$\vec{d}_{\Pi^1}$ is given by
$$
k=m_{i_1}+m_{i_2}+\cdots + m_{i_s}\,.
$$
\end{lemma}
\noindent {\bf Proof}. The depth $k$ of the gradation is equal to the maximal degree $d(\alpha)$,
$\alpha$ being a root. If $\alpha =k_1\alpha_1+\cdots +k_\ell\alpha_\ell$ is the
decomposition of a root $\alpha$ with respect to simple roots, then
$$
d(\alpha)=k_{i_1}+\cdots +k_{i_s}\leq d(\mu)=m_{i_1}+\cdots +m_{i_s}=k\,.
$$
$\qquad\Box$
\par\noindent
{\bf Irreducible submodules of the $\ggg^0$-module $\ggg^1$.}\ \ %%%%%%%%%%%%%%%%%%%%%%%%%%%%%%%%%%%%%%%%%%%%%%%%%%%%%%%%%%%%%%
%%%%%%%%%%%%%%%%%%%%%%%%%%%%%%%%%%%%%%%%%%%%%%%%%%%%%%%%%%%%%%%%%%%%%%%%%%%%%%%%%%%%%%%%%%%%%%%%%%%%%%%%%%%%%%%%%%%%
%%%%%%%%%%%%%%%%%%%%%%%%%%%%%%%%%%%%%%%%%%%%%%%%%%%%%%%%%%%%%%%%%%%%%%%%%%%%%%%%%%%%%%%%%%%%%%%%%%%%%%%%%%%%%%%%%%%%
Let $ \ggg=\sum\ggg^i $ be a fundamental gradation of a complex
semisimple Lie algebra $ \ggg $, defined by a label vector
$ \vec{d}$. Following \cite{GOV}, we describe the decomposition of
a $ \ggg^0 $-module into irreducible submodules. Set
$$
R^i=\{\alpha\in R\,\,\vert\,\, d(\alpha)=i\}=\{\alpha\in R\,\,\vert\,\,
\ggg_\alpha\subset\ggg^i\}
$$
and
$$
\Pi^i=\Pi\cap R^i =\{\alpha\in\Pi\,\,\vert\,\,d(\alpha)=i\}\,.
$$
For any simple root $ \gamma\in\Pi $, we put
$$
R(\gamma)=  \{\gamma  +(R^0 \cup \{0 \})\}\cap R  =
\{\alpha =\gamma  + \phi^0\in R,\,\,\,\phi^0\in R^0 \cup \{0 \} \}.
$$
We associate to any set of roots $ Q \subset R $ a subspace
$$
\ggg (Q)=\sum_{\alpha\in Q}\ggg_\alpha\subset \ggg \,.
$$
\begin{prop}{\rm (}\cite{GOV}{\rm )}\label{decomposizione}
The decomposition of a $\ggg^0$-module $\ggg^1$ into irreducible
submodules is given by
$$
\ggg^1=\sum_{\gamma\in\Pi^1}\ggg(R(\gamma))\,.
$$
Moreover, $\gamma$ is a lowest weight of the irreducible
submodule $\ggg(R(\gamma))$. In particular, the number of the irreducible
components is equal to the number $\#\Pi^1$ of the simple roots of degree
{\rm 1}.
\end{prop}
Since the $\ggg^0$-modules $\ggg^i$ and $\ggg^{-i}$ are dual,
Proposition \ref{decomposizione} gives also the decomposition of
the $\ggg^0$-module $\ggg^{-1}$ into irreducible submodules.
\section{Fundamental gradations of a real semisimple Lie algebra}%%%%%%%%%%%%%%%%%%%%%%%%%%%%%%%%%%%%%%%
%%%%%%%%%%%%%%%%%%%%%%%%%%%%%%%%%%%%%%%%%%%%%%%%%%%%%%%%%%%%%%%%%%%%%%
%%%%%%%%%%%%%%%%%%%%%%%%%%%%%%%%%%%%%%%%%%%%%%%%%%%%%%%%%%%%%%%%%%%%%%
\subsection{Real forms of a complex semisimple Lie algebra}
Now we recall the description of a real form of a
complex semisimple Lie algebra in terms of Satake diagrams. It is
sufficient to do this for complex simple Lie algebras.

Any real form of a complex semisimple Lie algebra $\ggg$ is the
fixed points set $\ggg^\sigma$ of an antilinear involution $\sigma$,
that is, an antilinear map $\sigma:\ggg\to \ggg$, which is an
automorphism of $\ggg$ as a real algebra, such that
$\sigma^2=\hbox{\rm id}$. We fix a Cartan decomposition
$$
\ggg^\sigma =\gk +\gm
$$
of the real form $\ggg^\sigma$, where $\gk$ is a maximal compact subalgebra of
$\ggg^\sigma$ and $\gm$ is its orthogonal complement
with respect to the Killing form $B$. Let
$$
\gh^\sigma =\gh_{\gk} +\gh_{\gm}
$$
be a Cartan subalgebra of $\ggg^\sigma$ which is consistent
with this decomposition and such that $\gh_{\gm}=\gh\cap\gm$
has maximal dimension. Then the root decomposition of
$\ggg^\sigma$, with respect to the subalgebra
$\gh^\sigma $, can be written as
$$
\ggg^\sigma =\gh^\sigma + \sum_{\lambda\in\Sigma}\ggg^\sigma_\lambda\,,
$$
where $\Sigma\subset (\gh^\sigma)^*$ is a (non-reduced) root system. The
number $m_\lambda =\dim\ggg_\lambda$ is the \emph{multiplicity}
of a root $\lambda\in\Sigma$.

Denote by $\gh =(\gh^\sigma)^\C$ the complexification of
$\gh^\sigma$ which is a $\sigma$-invariant Cartan subalgebra.
We denote by $\sigma^*$ the induced antilinear action of
$\sigma$ on $\gh^*$ given by
$$
\sigma^*\alpha =\overline{\alpha\circ\sigma}\,,\quad\alpha\in\gh^*\,.
$$
Consider the root space decomposition
$$
\ggg =\gh+\sum_{\alpha\in R} \ggg_\alpha
$$
of the Lie algebra $\ggg$ with respect to the Cartan subalgebra $\gh$. Note
that $\sigma^*$ preserves the root system $R$, i.e. $\sigma^*R=R$.
Now we relate the root space decomposition of $\ggg^\sigma$ and $\ggg$. We
define the subsystem of compact roots $R_\bullet$ by
$$
R_\bullet =\{\alpha\in R\,\,\vert\,\, \sigma^*\alpha = -\alpha \}=
\{\alpha\,\,\vert\,\,\alpha(\gh_\gm )=0\}
$$
and denote by $R'=R\setminus R_\bullet$ the complementary set of
non-compact roots. We can choose a system $\Pi$ of simple roots of
$R$ such that the corresponding system of positive roots $R_+$
satisfies the condition: $R'_+=R'\cap R_+$ is $\sigma$-invariant.
In this case, $\Pi$ is called a $\sigma$-\emph{fundamental system}
of roots. We denote by $\Pi_\bullet =\Pi\cap R_\bullet$ the set of
compact simple roots (which are also called black) and by $\Pi'
=\Pi\setminus \Pi_\bullet$ the non-compact simple roots (called
white). The action of $\sigma^*$ on white roots satisfies the
following property:

for any $\alpha\in\Pi'$ there exists a unique $\alpha'\in\Pi'$
such that\ $\sigma^*\alpha-\alpha'$\ is a linear combination of
black roots, i.e.
$$
\sigma^*\alpha = \alpha' +\sum_{\beta\in\Pi_\bullet}k_\beta\beta,
\quad k_\beta\in{\mathbb N}\,.
$$
In this case, we say that the roots $\alpha, \,\alpha'$ are
$\sigma$-\emph{equivalent} and we will write $\alpha\sim\alpha'$.
The information about fundamental  system ($\Pi =\Pi_\bullet \cup
\Pi'$) together with the $\sigma$-equivalence  can be visualized in
terms of the \emph{Satake diagram}\index{Satake diagram}, which is defined as follows.
\newline On the Dynkin diagram of the system of simple roots $\Pi$,
we paint the vertices which correspond to black roots into black and
we join the vertices which correspond to
$\sigma$-equivalent roots $\alpha,\,\alpha'$ by a curved arrow.\newline
By a slight abuse of notation, we will refer to the $\sigma$-fundamental system
$\Pi=\Pi_\bullet\cup\Pi'$, together with the $\sigma$-equivalence
$\sim$, as the \emph{Satake diagram}. This diagram is determined
by the real form $\ggg^\sigma$ of a complex simple Lie algebra $\ggg$
and does not depend on the choice of a Cartan subalgebra and a
$\sigma$-fundamental system. The list of Satake diagram of real
simple Lie algebras is known (see e.g. \cite{GOV}).\newline
Conversely, Satake diagram ($\Pi=\Pi_\bullet\cup\Pi',\sim$) allows
to reconstruct the action of $\sigma^*$ on $\Pi$, hence on $\gh^*$.
This action can be canonically extended to the antilinear involution
$\sigma$ of the complex Lie algebra $\ggg$. Hence,
{ \it there is a
natural $1-1$ correspondence between Satake diagrams subordinated to
the Dynkin diagram of a complex semisimple Lie algebra $\ggg$, up to
isomorphisms, and real forms $\ggg^\sigma$ of $\ggg$,
up to conjugations.} \newline
%We will describe real forms $\ggg^\sigma$ of a complex simple Lie algebra
% in terms of Satake diagrams ($\Pi=\Pi_\bullet\cup\Pi', \sim$).
\subsection{Gradations of a real semisimple Lie algebra}%%%%%%
%%%%%%%%%%%%%%%%%%%%%%%%%%%%%%%%%%%%%%%%%%%%%%%%%%%%%%%%%%%%%%%%%%%%%%%%%%%%%
%%%%%%%%%%%%%%%%%%%%%%%%%%%%%%%%%%%%%%%%%%%%%%%%%%%%%%%%%%%%%%%%%%%%%%%%%%%%%
Let $\ggg$ be a complex simple Lie algebra and $\ggg^\sigma$ be a real form of
$\ggg$ with a Satake diagram ($\Pi=\Pi_\bullet\cup\Pi',\sim)$.
\begin{comment}
We identify
$\Pi=\{\alpha_1,\ldots ,\alpha_\ell\}$ with  a $\sigma$-fundamental system,
 which is a system of simple roots of $\ggg$ with respect to a
Cartan subalgebra $\gh$ and  $\Pi_\bullet$ and $\Pi'$  with the
set of black and white roots respectively.\\
\end{comment}
Let $\vec{d}=(d_1,\ldots ,d_\ell)$ be a label vector of the simple roots
system $\Pi$ and $\ggg =\sum_{i\in\Z}\ggg^{i}$ be the corresponding
gradation of $\ggg$, with the grading element $d\in\gh\subset\ggg$.
\newline The following theorem gives necessary and sufficient
conditions in order that this gradation induces a gradation
$$
\ggg^\sigma =\sum_{i\in\Z}\ggg^\sigma\cap \ggg^i
$$
of the real form $\ggg^\sigma$. This means that
the grading element $d$ belongs to $\ggg^\sigma$. We denote by
$\Pi^0\subset\Pi$ the set of simple roots with label zero.
\begin{thm} $($\cite{D}$)$\label{djokovic}
A gradation of a complex semisimple Lie algebra $\ggg$,
associated with a label vector $\vec{d}=(d_1,\ldots ,d_\ell)$,
induces a gradation of the real form $\ggg^\sigma$, which
corresponds to a
Satake diagram $(\Pi = \Pi_{\bullet}\cup \Pi',\sim)$ if and only if
the following two conditions hold:
\begin{enumerate}
\item[i{\rm )}] $\Pi_\bullet \subset \Pi^0$, i.e. any black vertex
of the Satake diagram has label zero;
\item[ii{\rm )}] if
$\alpha\sim\alpha'$ for
$\alpha,\,\alpha'\in\Pi\setminus\Pi_\bullet$, then $d(\alpha)
=d(\alpha')$, i.e. white vertices of the Satake diagram which are
joint by a curved arrow have the same label.
\end{enumerate}
\end{thm}
\par
A label vector $\vec{d}=(d_1,\ldots ,d_\ell)$ of a Satake diagram
$(\Pi=\{\alpha_1,\ldots ,\alpha_\ell\}=\Pi_\bullet\cup\Pi',\sim)$
and the corresponding gradation of $\ggg$ are called of
{\em real type}
if they satisfy conditions i) and ii) of the theorem above, that is
black vertices have label zero and vertices related by a curved
arrow have the same label. Hence, we can state
Theorem \ref{djokovic} as follows
\begin{cor}
There exists a natural $1-1$ correspondence between label vectors
$\vec{d}$ of real type of a Satake diagram of a real semisimple Lie
algebra $\ggg^\sigma$ and gradations of $\ggg^\sigma$. The gradation
of $\ggg^\sigma$ is fundamental if and only if the corresponding
gradation of $\ggg$ is fundamental, i.e. $\vec{d} =
\vec{d_{\Pi^1}}$.
\end{cor}
\noindent
{\bf Irreducible submodules of the $\ggg^0$-module $\ggg^1$.}
Let $\ggg=\sum\ggg^i$ be a gradation of a complex
semisimple Lie algebra $\ggg$ with grading element $d$ and
$\ggg^\sigma =\sum(\ggg^\sigma)^i=\sum\ggg^i\cap\ggg^\sigma$ be a real form of
$\ggg$, consistent with this gradation. We denote by
$\,\,\,(\Pi=\Pi_\bullet
\cup\Pi',\,\sim)\,\,\,$ the Satake diagram of $\ggg^\sigma$.\newline
By Proposition \ref{decomposizione}, the decomposition of $\ggg^1$ into irreducible
$\ggg^0$-submodules is given by
$
\ggg^1=\sum_{\gamma\in\Pi^1}\ggg(R(\gamma))\,,
$
where $\Pi^1$ is the set of simple roots of label one. The
following obvious proposition describes the decomposition of
$(\ggg^\sigma)^0$-module $(\ggg^\sigma)^1$ into irreducible
submodules.
\begin{prop}
For any simple root $\gamma \in\Pi^1$ of label one,
there are two possibilities:
\begin{enumerate}
\item[i)] $\sigma^*\gamma=\gamma+\sum_{\beta\in\Pi_\bullet}k_\beta\beta$.
Then $\sigma^*\gamma \in R(\gamma)$ and the
$\ggg^0$-module $\ggg(R(\gamma))$ is $\sigma$-invariant;
\item[ii)] $\sigma^*\gamma =\gamma'+\sum_{\beta\in\Pi_\bullet}k_\beta\beta$,
where $\gamma\neq\gamma'\in \Pi^1$. Then, $\sigma^*R(\gamma)=R(\gamma')$ and
the two irreducible $\ggg^0$-modules
$\ggg(R(\gamma))$ and $\ggg(R(\gamma'))$ determine
one irreducible submodule
$
\ggg^\sigma\cap(\ggg(R(\gamma))+\ggg(R(\gamma')))
$
of $\ggg^\sigma$.
\end{enumerate}
\end{prop}

\begin{cor} \label{IrreducibleSubmodules}
Let $\ggg^\sigma =\sum(\ggg^\sigma)^i$ be the gradation of a real
semisimple Lie algebra  $\ggg^{\sigma}$, associated with a
label vector $\vec d $ of real type. Then irreducible
submodules of the $(\ggg^\sigma)^0$-module $(\ggg^\sigma)^{-1}$
correspond to vertices $\gamma$ with label one without curved
arrow and to pairs $(\gamma,\gamma')$ of vertices with label one
related by a curved arrow.
In particular, a decomposition of the $(\ggg^\sigma)^0$-module
$(\ggg^\sigma)^{-1}$  is determined by a decomposition
of the set $\Pi^1$ of vertices with label 1 into a disjoint union
$\Pi^1 = \Pi^1_+ \cup \Pi^1_- $  such that equivalent vertices belong
to the same component. The corresponding submodules $(\ggg^\sigma)^{-1}_+$
  and $(\ggg^\sigma)^1_-$ are  given  by
\begin{equation}\label{decomposition}
 (\ggg^\sigma)^{-1}_{\pm} = \ggg^\sigma \cap \sum_{\gamma\in \Pi^1_{\pm}} \ggg
  (R(-\gamma)).
  \end{equation}
\end{cor}
 We will always assume that a decomposition of $\Pi^1$ satisfies
 the above property.

\section{Classification of Maximally homogeneous para-CR manifolds} %(section 7)  %%%%%%%%%%%%%%%
%%%%%%%%%%%%%%%%%%%%%%%%%%%%%%%%%%%%%%%%%%%%%%%%%%%%%%%%%%%%%%%%%%%%%%%%%
%%%%%%%%%%%%%%%%%%%%%%%%%%%%%%%%%%%%%%%%%%%%%%%%%%%%%%%%%%%%%%%%%%%%%%%%%
 Let $\ggg^{\sigma}$ be a real semisimple Lie algebra associated with a
 Satake diagram $(\Pi = \Pi_{\bullet}\cup \Pi', \sim)$ with  the fundamental gradation
  defined by a subset $\Pi^1 \subset \Pi'$ and $F= \tilde G/\tilde P$ be the
  associated flag manifold.
  %Invariant almost para-CR structures on  $F=G/P$ correspond
  %to disjoint decompositions
  %$\Pi^1 = \Pi^1_+ \cup \Pi^1_-$ such that equivalent vertices
  %belong to the same component.

 By Theorem \ref{reductiontheorem}, an almost para-CR structure on $F = \tilde G/\tilde P$ associated with
 a decomposition  $\Pi^1 = \Pi^1_+ \cup \Pi^1_-$ is integrable (i.e. a para-CR structure)
 if and only if
 the $(\ggg^\sigma)^0$-submodules  $(\ggg^\sigma)^{-1}_+$ and $(\ggg^\sigma)^{-1}_-$ given by
 (\ref{decomposition}) are Abelian subalgebras
 of $\ggg^\sigma$. In order to give an integrability criterion, we introduce the following definitions.
\begin{definition}\label{admissible} Let $R$ be a system of roots and $\Pi$ be a system of simple roots. A subset
$\Pi^1 \subset \Pi$ is said to be {\rm admissible} if $\Pi^1$ contains at least two roots and
there are no roots of $R$ of the form
\begin{equation}\label{twoalphacondition}
2\alpha + \sum k_i \phi_i\,, \,\,\,\mbox{with}\,\, \alpha \in \Pi^1\,,\, \phi_i \in \Pi_0
= \Pi \setminus \Pi^1.
\end{equation}
\end{definition}

\begin{definition} Let $\ggg^\sigma$ be a  real semisimple Lie algebra
with a fundamental gradation defined by a subset $\Pi^1 \subset \Pi'$.
 We say that a decomposition $\Pi^1 = \Pi^1_+ \cup \Pi^1_-$
is {\em alternate} if the following conditions hold:
\begin{itemize}
\item[i)] if $\alpha\in\Pi^1_\pm$  and $\alpha'\sim\alpha$, then $\alpha'\in\Pi^1_\pm$;
\item[ii)] the vertices  in $\Pi^1_+$  and $\Pi^1_-$ appear  in the Satake diagram
 in alternate order. This means that each connected component of the graph obtained
 deleting vertices in $\Pi^1_+$ $($respectively in $\Pi^1_-$$)$ has not more than one vertex in
 $\Pi^1_-$ $($respectively in $\Pi^1_+$$)$.
\end{itemize}
\end{definition}
We are ready to state the following
\begin{prop} \label{alternateprop}Let $\ggg^\sigma$ be a semisimple real Lie algebra
with the fundamental gradation associated with a subset $\Pi^1 \subset \Pi$ and $F = \tilde G/\tilde P$
the associated flag manifold.
  A decomposition $\Pi^1 = \Pi^1_+ \cup \Pi^1_-$
defines a
para-CR structure on the flag manifold $F$ if and
only if  the subset $\Pi^1$ is admissible  and  the decomposition of $\Pi^1$ is alternate.
\end{prop}

 For the proof we need the following lemma.

\begin{lemma}\label{lemmaAbelCond} The subspace $\ggg^{1}_+ = \sum_{\gamma \in \Pi^1_+}\ggg(R(\gamma))$
$($hence also the subspace $(\ggg^\sigma)^{1}_+  = \ggg^\sigma \cap \ggg^1_+)$ which corresponds to a
subset $\Pi^1_+ \subset \Pi^1$ is an Abelian  subalgebra if and only if
there is no root $\beta$ of the form
\begin{equation}\label{AbelianCondition}
   \beta= \alpha + \alpha' + \sum k_i\phi_i
\end{equation}
where $\alpha, \alpha' \in \Pi^1_+$ and $\phi_i \in \Pi^0.$
The case $\alpha = \alpha'$ is allowed.
\end{lemma}
\noindent {\bf Proof.} If such a root $\beta$ exists, then $ [\ggg(R(\alpha), \ggg (R(\alpha'))] \neq 0$
and $\ggg^1_+$ is not an Abelian subalgebra. The converse is also clear.
$\Box$\smallskip\par\noindent
\noindent {\bf Proof of Proposition} \ref{alternateprop}.
  Let $ \Pi^1= \Pi^1_+ \cup \Pi^1_- $ be a decomposition of $\Pi^1$.  The condition
 (\ref{AbelianCondition}) for $\alpha = \alpha'$ is fulfilled if and only if
  $\Pi^1$ is admissible.  Assume now  that two different vertices
  $\alpha, \alpha'$ in $\Pi^1_+$ are
 connected in the Satake diagram by vertices in $\Pi^0=\Pi\setminus\Pi^1$. Then there is
 a root of the form (\ref{AbelianCondition}) and $\ggg^1_+$ is not a
 commutative  subalgebra. This shows that  the decomposition  which defines a para-CR structure
 on $F$ must be alternate.\\
 Conversely, assume that the decomposition is alternate. Then any two
 vertices $\alpha, \alpha' \in \Pi^1_+$ belong to different
 connected components of the  Satake graph with deleting
 $\Pi^1_-$. This implies that  there is no root of the form
 (\ref{AbelianCondition}) for $\alpha \neq \alpha'$. Then
  Lemma \ref{lemmaAbelCond} shows that $(\ggg^\sigma)^1_+$ is a commutative
  subalgebra.
 The same argument is applied also for $(\ggg^\sigma)^1_-$.
 $\Box$\smallskip\par\noindent
%%%%%%%%%%%%%%%%%%%%%%%%%%%%%%%%%%%%%%%%%%%%%
\begin{comment}
The following Proposition describes  admissible subsystems $\Pi^1$
 of a system $\Pi$ of simple roots for any  indecomposable root
 system $R$.\\
\end{comment}
We enumerate simple roots of complex simple Lie $\ggg$
algebras as in \cite{Bou}. Let $\Pi =\{\alpha_1, \ldots , \alpha_{\ell}\}$ be
the simple roots of $\ggg$, which are
identified with vertices of the corresponding Dynkin diagram. We denote the elements of
a subset $\Pi^1 \subset \Pi$ (respectively $\Pi^1 \subset \Pi'$)  which defines a fundamental gradation
of $\ggg$ (respectively $\ggg^\sigma$) by
$$
\alpha_{i_1}, \ldots , \alpha_{i_k},\,\,\,\,\, i_1 < i_2 < \cdots < i_k.
$$
\begin{prop}\label{lista-complessa} Let $\Pi$ be a system of simple roots of a root system $R$ of a complex
simple Lie algebra $\ggg$. Then
a subset $\Pi^1 \subset \Pi$ of at least two elements is admissible $($see Definition \ref{admissible}$)$ in the following cases:

\begin{itemize}
\item  for $ \ggg = A_\ell$, in all cases;
\item  for $\ggg = B_\ell$, under the
condition: $i_{k}=i_{k-1}+1$;
\item  for $\ggg = C_\ell$, under the
condition: $i_{k}=\ell$;
\item for $\ggg = D_\ell$, under the condition: if $i_k < \ell-1$,
then $i_k=i_{k-1}+1$;
\item  for $\ggg = E_6$, in all cases except the following ones:
$$
\{\alpha_1,\alpha_4\}\,,\,\{\alpha_1,\alpha_5\}\,,\,\{\alpha_3,\alpha_6\}\,,\,\{\alpha_4,\alpha_6\}\,,\,
\{\alpha_1,\alpha_4,\alpha_6\}\,;
$$
\item for $\ggg = E_7$, in all cases except the following ones:
\begin{eqnarray*}
&\{\alpha_1,\alpha_4\}\,,\,\{\alpha_1,\alpha_5\}\,,\,\{\alpha_3,\alpha_6\}\,,\,
\{\alpha_4,\alpha_6\}\,,\,\{\alpha_1,\alpha_6\}\,,\,&\\
&\{\alpha_2,\alpha_7\}\,,\,\{\alpha_3,\alpha_7\}\,,\,\{\alpha_4,\alpha_7\}\,,\,
\{\alpha_5,\alpha_7\}\,,\,&\\
&\{\alpha_1,\alpha_4,\alpha_6\}\,,\,\{\alpha_1,\alpha_4,\alpha_7\}\,,\,\{\alpha_1,\alpha_5,\alpha_7\}\,,\,
\{\alpha_3,\alpha_6,\alpha_7\}\,,\,\{\alpha_4,\alpha_6,\alpha_7\}\,,\,&\\
&\{\alpha_1,\alpha_4,\alpha_6,\alpha_7\}\,;&
\end{eqnarray*}
\item for  $\ggg = E_8$, in all cases except the following ones:
\begin{eqnarray*}
&\{\alpha_1,\alpha_4\}\,,\,\{\alpha_1,\alpha_5\}\,,\,\{\alpha_3,\alpha_6\}\,,\,
\{\alpha_4,\alpha_6\}\,,\,\{\alpha_1,\alpha_6\}\,,\,&\\
&\{\alpha_2,\alpha_7\}\,,\,\{\alpha_3,\alpha_7\}\,,\,\{\alpha_4,\alpha_7\}\,,\,
\{\alpha_5,\alpha_7\}\,,\,&\\
&\{\alpha_1,\alpha_4,\alpha_6\}\,,\,\{\alpha_1,\alpha_4,\alpha_7\}\,,\,\{\alpha_1,\alpha_5,\alpha_7\}\,,\,
\{\alpha_3,\alpha_6,\alpha_7\}\,,\,\{\alpha_4,\alpha_6,\alpha_7\}\,,\,&\\
&\{\alpha_1,\alpha_4,\alpha_6,\alpha_7\}\,,\,&\\
&\{\alpha_1,\alpha_7\}\,,\,\{\alpha_1,\alpha_8\}\,,\,\{\alpha_2,\alpha_8\}\,,\,\{\alpha_3,\alpha_8\}\,,\,\{\alpha_4,\alpha_8\}
\,,\,\{\alpha_5,\alpha_8\}\,,\,\{\alpha_6,\alpha_8\}\,,&\\
&\{\alpha_1,\alpha_4,\alpha_8\}\,,\,\{\alpha_1,\alpha_5,\alpha_8\}\,,\,\{\alpha_3,\alpha_6,\alpha_8\}\,,\,
\{\alpha_4,\alpha_6,\alpha_8\}\,,\,\{\alpha_1,\alpha_6,\alpha_8\}\,,\,&\\
&\{\alpha_2,\alpha_7,\alpha_8\}\,,\,\{\alpha_3,\alpha_7,\alpha_8\}\,,\,\{\alpha_4,\alpha_7,\alpha_8\}\,,\,
\{\alpha_5,\alpha_7,\alpha_8\}\,,\,&\\
&\{\alpha_1,\alpha_4,\alpha_6,\alpha_8\}\,,\,\{\alpha_1,\alpha_4,\alpha_7,\alpha_8\}\,,\,
\{\alpha_1,\alpha_5,\alpha_7,\alpha_8\}\,,\,
\{\alpha_3,\alpha_6,\alpha_7,\alpha_8\}\,,&\\
&\{\alpha_4,\alpha_6,\alpha_7,\alpha_8\}\,,\,&\\
&\{\alpha_1,\alpha_4,\alpha_6,\alpha_7,\alpha_8\}\,;&
\end{eqnarray*}
\item for $ \ggg = F_4 $, in all cases except the following ones:
$$
\{\alpha_1,\alpha_3\}\,,\,\{\alpha_1,\alpha_4\}\,,\,\{\alpha_2,\alpha_4\}\,,\,\{\alpha_3,\alpha_4\}\,,\,
\{\alpha_1,\alpha_3,\alpha_4\}\,;
$$
\item for $ \ggg = G_2$, in the  case $\{\alpha_1,\alpha_2\}\,.$
\end{itemize}
In cases different from $D_\ell$, $E_6$, $E_7$ and $E_8$, for any $\Pi^1$
given as above it is possible to give an alternate decomposition $\Pi^1=\Pi^1_+\cup\Pi^1_-$.\newline
For $D_\ell$, an alternate decomposition of $\Pi^1$ can be given in the following cases:
\begin{itemize}
\item $\alpha_{\ell-2}\in \Pi^1$,
\item $\Pi^1$ is contained in at most two of the branches issuing from $\alpha_{\ell-2}$.
\end{itemize}
For $E_6$, $E_7$ and $E_8$, an alternate decomposition of $\Pi^1$ can be given in the following cases:
\begin{itemize}
\item $\alpha_{4}\in \Pi^1$,
\item $\Pi^1$ is contained in at most two of the branches issuing from $\alpha_{4}$.
\end{itemize}
\end{prop}
\noindent {\bf Proof.} We have to describe all subsets $\Pi^1$
  of $\Pi$ which  satisfy
  (\ref{twoalphacondition}). This condition can be
  reformulated as  follows. For any $\alpha \in \Pi^1$, denote by $\Pi_{\alpha}$ the
  connected component of the subdiagram of
   the Dynkin diagram $\Pi$ obtained by deleting vertices in $\Pi^1 \setminus
  \{\alpha\}$ and containing $\alpha$. Then the root system associated with  $\Pi_\alpha$
  has no roots of the form
  $$ \beta = 2\alpha + \sum_{\phi \in \Pi_\alpha \setminus \{\alpha\}}k_{\phi}\phi.$$
 Using this condition and the decomposition of any root
 into a linear combination of simple roots, one can prove the
 proposition. \\
 In the case of $A_\ell$, any root has coefficient $0,1$ in the
 decomposition into simple roots. Hence, any
 decomposition satisfies the property (\ref{twoalphacondition}).\par\smallskip\noindent
 In the case of $B_\ell$, any root which has coefficient 2 has the
 form
$$
\sum_{i\leq h< j}\alpha_h +2\sum_{j\leq h\leq\ell}\alpha_h\,,\qquad(1\leq i<j\leq\ell)\,.
$$
Hence the condition (\ref{twoalphacondition}) holds if and only
if the last two roots in $\Pi^1$ are consecutive, i.e. $i_{k-1}+1=i_k$.\par\smallskip\noindent
In the case of $C_\ell$, the roots  with a coefficient 2 are given by
\begin{eqnarray*}
&{}&\sum_{i\leq h<j}\alpha_h +2\sum_{j\leq h<\ell}\alpha_h+\alpha_\ell\,,\qquad(1\leq i<j\leq\ell)\,,\\
&{}&2\sum_{i\leq h<\ell}\alpha_h +\alpha_\ell\,,\qquad(1\leq i<\ell)\,.
\end{eqnarray*}
The second formula implies that there are no roots of the form given in \eqref{twoalphacondition}
if and only if $i_k=\ell$.\par\smallskip\noindent
 In the case of $D_\ell$, the roots with a coefficient $2$ are
$$
\sum_{i\leq h<j}\alpha_h +2\sum_{j\leq h<\ell-1}\alpha_h+\alpha_{\ell-1}+\alpha_\ell\,,\qquad(1\leq
i<j<\ell-1)\,.
$$
The condition (\ref{twoalphacondition}) fails if and only if
the last two roots $\alpha_{i_{k-1}}\,, \alpha_{i_k}$ satisfy $i_{k-1}< i_k -1$ and $i_k <
\ell-1$.\par\smallskip\noindent
 The case of exceptional Lie algebras  can be  treated in a similar way, by using tables in \cite{Bou}.
$\Box$\medskip\par\noindent
Let $\Pi^1\subset\Pi'$ be an admissible subset which defines a fundamental gradation of $\ggg^\sigma$.
An alternate decomposition of $\Pi^1=\Pi^1_+\cup\Pi^1_+$ can be given if the conditions of
Proposition \ref{lista-complessa} are satisfied and, in addition, the following ones hold:
\begin{itemize}
\item for $\su(p,q)$, it has to be $q=p$ and $\alpha_p\in\Pi^1$;\medskip
\item for $\so(\ell-1,\ell+1)$, it has to be $\Pi^1\cap\{\alpha_{\ell -1},\alpha_\ell\}=\emptyset$ or
$\{\alpha_{\ell -2},\alpha_{\ell -1},\alpha_\ell\}\subset\Pi^1$;\medskip
%\item for $\so^*(2\ell)$, it has to be $\Pi^1\cap\{\alpha_{\ell -1},\alpha_\ell\}=\emptyset$;
\item for $E_6$II, it has to be $\alpha_4\in\Pi^1$ and if $\alpha_2\notin\Pi^1$, then
$\{\alpha_3,\alpha_5\}\subset\Pi^1$;
\end{itemize}
while for $\so^*(2\ell)$ and $E_6$III there is no alternate decomposition of $\Pi^1$.
\medskip\par\noindent
Proposition \ref{alternateprop} implies the following final theorem.
\begin{thm} Let $(\Pi = \Pi_{\bullet} \cup \Pi', \sim)$ be a
Satake diagram of a simple real Lie algebra $\ggg^\sigma $ and $\Pi^1 \subset \Pi'$
be an admissible subset as described above.
Let $\tilde G$ be the simply connected Lie group
with the Lie algebra $\ggg^\sigma$ and $\tilde P$ be the parabolic
subgroup of $\tilde G$ generated by the non-negatively graded subalgebra
$$
\gp = \sum_{i \geq 0} (\ggg^\sigma)^i
$$
associated with the grading element $\vec d_{\Pi^1}$. Then the alternate decomposition
$\Pi^1 = \Pi^1_+ \cup \Pi^1_-$ defines a decomposition
$$(\ggg^\sigma)^1 = (\ggg^\sigma)^1_+ +  (\ggg^\sigma)^1_- $$
of the $(\ggg^\sigma)^0$-module $(\ggg^\sigma)^1$ into a sum of two
commutative subalgebras. This decomposition determines an
invariant para-CR structure on the
simply connected  flag manifold $F = \tilde G/\tilde P$. Moreover,  any simply
connected maximally homogeneous para-CR manifolds of semisimple
type is a direct product of such manifolds.
\end{thm}
\par
\noindent {\bf Acknowledgement.} The authors would like to thank P. Nurowski for bringing their attention to the paper \cite{NS}
and for useful discussions.

\printindex
\end{document}